\documentclass[12pt]{article}

\usepackage{graphicx}
\usepackage{amssymb}
\usepackage{amsfonts}

\begin{document}

\title{Contractions and generalized Casimir invariants}

\author{Rutwig Campoamor-Stursberg\\ Depto. Geometr\'{\i}a y Topolog\'{\i}a\\Fac. CC. Matem\'aticas U.C.M.\\
E-28040 Madrid ( Spain)\\e-mail: rutwig@nfssrv.mat.ucm.es}

\date{}

\maketitle

\begin{abstract}
We prove that if $\frak{g}^{\prime}$ is a contraction of a Lie algebra $\frak{g}$ then the number of functionally independent invariants of $\frak{g}^{\prime}$ is at least that of $\frak{g}$. This allows to determine explicitly the number of invariants of Lie algebras carrying a supplementary structure, such as linear contact or linear forms whose differential is symplectic.   

\bigskip

\textit{Keywords:} Invariants, contractions, Lie algebras\newline 
PACS bumber: 02.20S
\end{abstract}

Invariants of the coadjoint representation, also called generalized Casimir invariants, are well known to characterize irreducible representations of Lie algebras, which justifies their importance in Physics, where invariants are used to characterize specific properties of a physical system, like mass formulas for dynamical groups. For the classical Lie algebras the problem was solved long ago by Racah, and by the Levi decomposition the class which must be analyzed is that of solvable Lie algebras. Various authors have approached the problem in recent years, and the difficulty of finding or even characterizing the invariants of solvable Lie algebras has been pointed out. Among the various interesting questions referring to these algebras, a characterization of solvable Lie algebras with abelian nilradical admitting only trivial invariants was presented in\cite{Nd}. In the same paper, the author commented the importance of finding a corresponding characterization 
for solvable Lie algebras with non-abelian nilradical. In view of the cases treated and other examples, such a characterization probably does not exist. However, sufficiency conditions to ensure the non-existence of non-trivial invariants can be found. The method we propose here is to study the problem considering contractions of Lie algebras, and to find a formula which expresses a relation between the number of functionally independent invariants of a Lie algebra and its contractions. This will be of much interest if we suppose that the contraction carries an additional structure, such as linear contact forms or symplectic forms. In particular, this will enable us to determine the number of invariants without any information about the precise structure of the contracting algebra.

\bigskip

An $n$-dimensional Lie algebra $\frak{g}$ may be considered as an element $\mu\in Hom\left(\bigwedge^{2}V,V\right)$, where $V$ is an $n$-dimensional vector space. The set $\mathcal{L}^{n}$ of Lie algebras is then a subset of the variety $Hom\left(\bigwedge^{2}V,V\right)$ on which the general linear group $GL\left(n,K\right)$ acts by changes of basis, i.e.,
\begin{eqnarray*}
(g\circ\mu)(x,y)=g^{-1}\left(\mu\left(gx,gy\right)\right),\quad g\in GL\left(n,K\right); x,y\in V
\end{eqnarray*}
Clearly the orbit of this action are the isomorphism classes of $\mu$. Now a Lie algebra $\mu_{1}$ is called a contraction of a Lie algebra $\mu_{0}$ if $\mu_{1}\in\overline{\mathcal{O}\left(\mu_{0}\right)}$, the Zariski closure of the orbit, Classically \cite{IW} this is expressed  by $\mu_{1}=\lim_{\epsilon\rightarrow 0}\left(g(\epsilon)\circ\mu_{0}\right)$. Therefore, if $\left\{C_{ij}^{k}\right\}$ are the structure constants of $\frak{g}_{0}=\left(\mathbb{K}^{n},\mu_{0}\right)$ over a basis $\left\{X_{1},..,X_{n}\right\}$ and $\left\{C_{ij}^{k}(\epsilon)\right\}$ the constants on the transformed basis by $g(\epsilon)$, the law of $\mu_{1}$ is given by 
\begin{eqnarray}
\widetilde{C_{ij}^{k}}=\lim_{\epsilon\rightarrow 0}C_{ij}^{k}(\epsilon)
\end{eqnarray}
As follows from the work of Beltrametti and Blasi \cite{Be}, the number of functionally independent invariants of the coadjoint representation $ad^{*}$ of a Lie algebra $\frak{g}$ is given by $\mathcal{N}=dim\left(\frak{g}\right)-r\left(\frak{g}\right)$, where 
$r\left(\frak{g}\right)$ is the maximum rank of the commutator table considered as a $\left(n\times n\right)$-matrix. The result we prove in this note is the following:
\begin{eqnarray}
\textrm{if}\quad \mu_{1}\quad \textrm{is a contraction of}\quad \mu_{0},\quad \textrm{then}\quad \mathcal{N}\left(\mu_{1}\right)\geq \mathcal{N}\left(\mu_{0}\right) 
\end{eqnarray}
It is clear that, starting from an arbitrary basis $\left\{X_{1},..,X_{n}\right\}$ of $\frak{g}_{0}$ we have:  
\begin{eqnarray}
rank\left(C_{ij}^{k}(\epsilon)\right)\geq rank\left(\lim_{\epsilon\rightarrow 0}\left(C_{ij}^{k}(\epsilon)\right)\right)
\end{eqnarray}
This equation holds for any basis and any $g(\epsilon)\in GL\left(n,\mathbb{K}\right)$. Now, as any contraction can be realized as a deformation \cite{Fi}, we obtain that 
\begin{eqnarray}
r\left(\mu_{0}\right)\geq r\left(\mu_{1}\right)
\end{eqnarray}
and from the formula for the number of invariants:
\begin{eqnarray}
\mathcal{N}\left(\mu_{0}\right)\leq \mathcal{N}\left(\mu_{1}\right)
\end{eqnarray}
This results coincides with the intuition that contractions have "less brackets" than the Lie algebra they come from. Geometrically this is clear, as the dimension of the orbit of contracted algebras is lower than the orbit dimension of the starting algebra, and therefore one should expect that the contraction has more invariants. Observe further that this result cannot be formulated in terms of deformations, since there exist deformations which are not related to a contraction \cite{Fi}. The following assertion are consequences of formula $(2)$:
\begin{enumerate}

\item if $\frak{g}=\left(\mathbb{K}^{2n+1},\mu\right)$ has a linear contact form, then $\mathcal{N}(\mu)\leq 1$.

\item if $\frak{g}$ is a frobeniusian Lie algebra, then has no non-trivial invariants.

\item if $\frak{g}$ has a contraction without no non-trivial invariants, then $\frak{g}$ itself has only trivial invariants.
\end{enumerate}
  
The first assertion follows immediately from the fact that any Lie algebra admitting a linear contact form contracts on a Heisenberg Lie algebra $\frak{h}_{n}$, which implies that the dimension of $\frak{g}$ is odd \cite{G2}. Now Heisenberg Lie algebras admit only one non-trivial invariant. For the second, frobeniusian Lie algebras, i.e., Lie algebras admitting a linear form whose differential is symplectic, are known (over the real or complex field) to contract on rank one solvable Lie algebras whose nilradical is isomorphic to a Heisenberg Lie algebra, called frobeniusian model Lie algebras \cite{G1}. 
A direct computation shows that frobeniusian model Lie algebras do not have non-trivial invariants. Observe that these model algebras constitute a particular case of the Lie algebras studied in \cite{Wi}, and that their number of functionally independent invariants follows easily by evaluation in the rank dependent formulae for $\mathcal{N}$ obtained in this reference.
By (iii), which is a trivial consequence of $(2)$, the result follows. \newline We also note the interest of this result for the study of contraction trees of Lie algebras, as equation $(2)$ establishes a necessary condition for an algebra to be a contraction.


\begin{thebibliography}{9}

\bibitem{Be} Beltrametti E G and Blasi A 1966 \textit{Phys. Lett.} \textbf{20} 62

\bibitem{Fi} Fialowski A and O'Halloran J 1988 \textit{Comm. Algebra} \textbf{112} 315

\bibitem{IW} In\"{o}n\"{u} E and Wigner E P 1953 \textit{Proc. Natl. Acad. Sci. U. S.} \textbf{39} 510

\bibitem{G1} Goze M 1981 \textit{C. R. A. S. Paris} \textbf{293} 425

\bibitem{G2} Goze M 1981 \textit{C. R. A. S. Paris} \textbf{292} 813

\bibitem{Nd} Ndogmo J P 2000 \textit{J. Phys. A: Math. Gen.} \textbf{33} 2273

\bibitem{Wi} Rubin J L and Winternitz P 1993 \textit{J. Phys. A: Math. Gen.} \textbf{26} 1123

\end{thebibliography}
\end{document}